\documentclass[11pt, letterpaper]{amsart}
\usepackage{amsmath}
\usepackage{amssymb}
\usepackage[left=1in,right=1in,bottom=0.8in,top=0.8in]{geometry}
\usepackage{amsfonts}
\usepackage{graphicx}
\usepackage[OT2,T1]{fontenc}
\usepackage[font=small,labelfont=bf]{caption}
\usepackage{epstopdf}
\usepackage{hyperref}
\usepackage[dvipsnames]{xcolor}
\usepackage{amsthm}
\usepackage{float}
\usepackage{pgfplots}
\usepackage{listings}
\usepackage{longtable}
\usepackage{mathrsfs}
\usepackage{bbm}
\usepackage{euscript}
\usepackage{mathtools}
\usepackage{tikz}
\usepackage{comment}
\usepackage{tikz-cd}
\usepackage{amsxtra}
\usepackage{pifont}
\usepackage{amsbsy}
\usepackage{wasysym}
\usepackage{calligra}
\usepackage{url}
\usepackage{enumerate}
\usepackage[
backend=biber,
style=alphabetic,
]{biblatex}

\makeatletter
\renewcommand{\tocsection}[3]{%
  \indentlabel{\@ifnotempty{#2}{\bfseries\ignorespaces#1 #2.\quad}}\bfseries#3}
  \renewcommand{\tocsubsection}[3]{%
  \indentlabel{\hspace{2.5em}\@ifnotempty{#2}{\ignorespaces#1 #2.\quad}}#3}
\makeatother

\newtheorem{theorem}{Theorem}[section]
\newtheorem{corollary}[theorem]{Corollary}
\newtheorem{proposition}[theorem]{Proposition}
\newtheorem{lemma}[theorem]{Lemma}
\newtheorem{conjecture}[theorem]{Conjecture}

\theoremstyle{definition}
\newtheorem{definition}[theorem]{Definition}

\theoremstyle{remark}
\newtheorem{remark}[theorem]{Remark}

\DeclareSymbolFont{cyrletters}{OT2}{wncyr}{m}{n}
\DeclareMathSymbol{\Sha}{\mathalpha}{cyrletters}{"58}

\definecolor{energy}{RGB}{114,0,172}
\definecolor{refblue}{RGB}{30,90,180}
\definecolor{freq}{RGB}{45,177,93}
\definecolor{spin}{RGB}{251,0,29}
\definecolor{signal}{RGB}{203,23,206}
\definecolor{circle}{RGB}{217,86,16}
\definecolor{average}{RGB}{203,23,206}

\colorlet{shadecolor}{gray!20}
\pgfplotsset{compat=1.9}

\usepgflibrary{fpu}

\newcommand{\RR}{\mathbb{R}}
\newcommand{\ZZ}{\mathbb{Z}}

\newcommand{\QQ}{\mathbb{Q}}
\newcommand{\CC}{\mathbb{C}}

\newcommand{\PP}{\mathbb{P}}

\newcommand{\FF}{\mathbb{F}}

\newcommand{\cE}{\mathcal{E}}

\newcommand{\cL}{\mathcal{L}}

\newcommand{\Aut}{\operatorname{Aut}}

\newcommand{\divs}{\operatorname{div}}

\newcommand{\cEnd}{\mathscr{E}\text{\kern -3pt {\calligra\large nd}}\,}
\newcommand{\cBun}{\mathscr{B}\text{\kern -3pt {\calligra\large un}}\,}
\newcommand{\cCoh}{\mathscr{C}\text{\kern -3pt {\calligra\large oh}}\,}

\newcommand{\cHom}{\mathscr{H}\text{\kern -3pt {\calligra\large om}}\,}

\newcommand{\ind}{\operatorname{ind}}

\newcommand{\NS}{\operatorname{NS}}

\newcommand{\pr}{\operatorname{pr}}

\newcommand{\rank}{\operatorname{rank}}

\hypersetup{
    colorlinks=true,
    linkcolor=refblue,
    citecolor=red,
}

\author[Ting Gong]{Ting Gong}
\email{tgong2@uw.edu}

\title{On potential density of rational points on hyperkahler varieties}
\begin{document}
                                                           
\begin{abstract}
In this paper, we prove that the rational points on hyperk\"ahler varieties in each of the four known deformation families admitting a rational Lagrangian fibration are potentially dense. Furthermore, we give several applications and construct examples of hyperk\"ahler varieties beyond these families satisfying the potential density property. Our approach is first to prove that rational points on hyperk\"ahler varieties admitting two rational Lagrangian fibrations with distinct associated isotropic rays are potentially dense. We then construct a projective Shafarevich--Tate twist admitting two such fibrations, together with a dominant rational map from the twist to the original variety.
\end{abstract}

\maketitle
\tableofcontents

\section{Introduction}
The study of the distribution of rational points has a long history. Two notable classical results are the Mordell--Weil theorem \cite{Weil_29} on the structure of rational points on abelian varieties and Faltings' theorem \cite{Faltings_83} on the finiteness of rational points on curves of genus at least two. By Faltings' theorem, such curves provide examples of varieties whose rational points are not Zariski dense over any finite extension of a number field. However, abelian varieties, including elliptic curves, behave differently. By \cite[Proposition~3.1]{Hassett_Tschinkel_00}, every abelian variety over a number field admits a nondegenerate algebraic point, whose cyclic subgroup is Zariski dense. After a finite extension making this point rational, we obtain a Zariski dense set of rational points. 

Let $X$ be a variety over a number field $k$. We say that rational points on $X$ is \emph{potentially dense} if there exists a finite field extension $L/k$ such that $X(L)$ is Zariski dense in $X_L$. This definition generalizes the preceding observations and emphasizes that potential density is an arithmetic property of $X$. Therefore, a natural question is: under what conditions are the rational points on varieties potentially dense?

One of the initial conjectures in this direction is the Bombieri--Lang conjecture \cite{Bombieri_90}, which generalizes Faltings' theorem and predicts that rational points on varieties of general type are not Zariski dense. In the opposite direction, Campana \cite{Campana_04} conjectured that rational points on special varieties are potentially dense. In his terminology, a special variety admits no fibration whose orbifold base is of general type. In particular, known examples of varieties satisfying the potential density property, such as abelian varieties and elliptic K3 surfaces \cite{Hassett_Tschinkel_00,Bogomolov_Tschinkel_00,Lai_Nakahara_22}, rational varieties such as $\PP^n$, and certain Calabi--Yau threefolds \cite{Bogomolov_Halle_Halvard_Tanimoto_18}, all belong to this category.

One of the main open cases is the potential density of rational points on K3 surfaces that do not admit genus-one fibrations, particularly those of geometric Picard rank one. Several works have studied this problem and related questions. Hassett and Tschinkel \cite{Hassett_Tschinkel_08} constructed examples of K3 surfaces of geometric Picard rank one with Zariski dense rational points over complex function fields, while van Luijk \cite{vanLuijk_07} constructed K3 surfaces of geometric Picard rank one over number fields with infinitely many rational points. More recently, Mart\'inez-Mar\'in \cite{MartinezMarin_26} established a uniform bound on the degree of a field extension needed to obtain infinitely many rational points on K3 surfaces of degree two. Furthermore, Bartsch, Javanpeykar, and Levin \cite{Bartsch_Javanpeykar_Levin_25} studied Campana's conjecture and established potential density results for symmetric powers of products of curves.

In this paper, we consider a natural higher-dimensional generalization of the potential density conjecture for K3 surfaces, which is also a special case of Campana's conjecture.

\begin{conjecture}\label{conj:hyperk}
    Every projective irreducible hyperk\"ahler variety over a number field has the potential density property.
\end{conjecture}

In particular, we study the following conjecture extensively, generalizing the potential density result for elliptic K3 surfaces of \cite{Bogomolov_Tschinkel_00}.

\begin{conjecture}\label{conj:one-fibration}
    Let $X$ be a projective irreducible hyperk\"ahler variety of dimension $2n$ over a number field $k$ admitting a rational Lagrangian fibration. Then the rational points on $X$ is potentially dense.
\end{conjecture}

In fact, in this paper, we establish Conjecture~\ref{conj:one-fibration} in almost all cases. Before presenting our results, we recall some basic definitions. A smooth projective variety $X$ over $\CC$ is called an \emph{irreducible hyperk\"ahler variety} if $X$ is simply connected and $H^0(X,\Omega_X^2)$ is generated by a nowhere degenerate holomorphic $2$-form \cite{Beauville_83}. One of the most notable open problems in this area is to obtain a complete classification of the deformation types of hyperk\"ahler varieties. To our knowledge, four types of deformation families are currently known: $K3^{[n]}$-type and generalized Kummer type constructed by Beauville \cite[Sections~6--7]{Beauville_83}, and the two exceptional deformation types in dimensions six and ten constructed by O'Grady \cite{OGrady_03,OGrady_99}, denoted by OG6 and OG10, respectively.

By a \emph{hyperk\"ahler variety} $X$ of dimension $2n$ over a number field $k$, we mean a smooth projective variety such that $X_\CC$ is an irreducible hyperk\"ahler variety. Furthermore, by a \emph{regular Lagrangian fibration} on $X$ over $k$, we mean a surjective morphism $f:X\to B$, where $B$ is a normal projective variety of dimension $n$, the fibres are geometrically connected, and the geometric generic fibre is Lagrangian \cite{Huybrechts_Mauri_22}. In particular, the smooth geometric fibres are abelian varieties. When $B$ is smooth, Hwang's theorem \cite{Hwang_08} implies that $B_\CC\simeq\PP^n_\CC$.

Our main result is the following theorem.

\begin{theorem}[=Theorem~\ref{thm:main}]
    Let $X$ be an irreducible hyperk\"ahler variety over a number field $k$ of $K3^{[n]}$-type or generalized Kummer type with $n\geq 2$, or of OG6 or OG10 type. If, after a finite extension of $k$, $X$ admits a rational Lagrangian fibration, then the rational points on $X$ is potentially dense.
\end{theorem}

We also prove the following cases: if $X$ has no geometric prime exceptional divisors and either $\rho(X_{\bar{k}})\geq 3$ or $\rho(X_{\bar{k}})=2$ with $b_2(X)\geq 5$, then the existence of a rational Lagrangian fibration implies that rational points on $X$ is potentially dense (=Theorem ~\ref{thm:fibration_rank_3}, ~\ref{thm:fibration_rank_2}). Furthermore, if $X$ admits a regular Lagrangian fibration with geometrically simple generic fibre and rational points on every geometric prime exceptional divisor of $X$ is potentially dense, then rational points on $X$ is potentially dense under the same Picard rank conditions (=Theorem ~\ref{thm:inductive_theorem}). We also construct examples of hyperk\"ahler varieties satisfying the potential density property in arbitrary dimensions that admit dominant generically finite rational maps from products of K3 surfaces and abelian varieties (=Theorem ~\ref{thm:example}).

We further deduce several applications to Conjecture~\ref{conj:hyperk}. In particular, we prove that rational points on $X$ is potentially dense if $\rho(X_{\bar{k}})\geq 5$ (=Corollary ~\ref{cor:Pic-5}), or if $X$ is of $K3^{[n]}$-type or generalized Kummer type and its extended algebraic lattice has Witt index $2$ (=Corollary ~\ref{thm:witt-two-density}).

\subsection{Methodology} 
We prove the above results as follows. We first show that rational points on any hyperk\"ahler variety admitting two rational Lagrangian fibrations with distinct associated isotropic rays is potentially dense. The idea essentially follows Hassett--Tschinkel \cite{Hassett_Tschinkel_00}. We first construct a dominant and generically finite rational map $X\dashrightarrow B_1\times B_2$, where $B_1$ and $B_2$ are the bases of the two Lagrangian fibrations. Then, considering a smooth fibre of the second fibration, which is geometrically an abelian variety, we obtain a potentially dense multisection of the first fibration. By choosing a nondegenerate point in the acting abelian variety of the first fibration, we obtain a rational section of the abelian scheme pulled back to this multisection. We then show that this section remains nondegenerate at the geometric generic point, which allows us to deduce potential density.

When $X$ has no geometric prime exceptional divisors and $\rho(X_{\bar{k}})\geq 3$, the Beauville--Bogomolov form on $\NS(X_{\bar{k}})_\RR$ has signature $(1,\rho(X_{\bar{k}})-1)$. By the action of the birational automorphism group and the Borel density theorem, we can find a movable isotropic ray distinct from the ray associated with the original fibration. When $\rho(X_{\bar{k}})=2$ and $b_2(X)\geq 5$, by the boundedness of MBM classes \cite{Amerik_Verbitsky_20} and the description of the K\"ahler cone \cite{Amerik_Verbitsky_15}, we show that the other isotropic boundary ray becomes nef on a smooth birational model. By the hyperk\"ahler SYZ theorem \cite{Engel_Mauri_26}, this gives another rational Lagrangian fibration, and we can apply our two-fibration theorem. If $X$ has a geometric prime exceptional divisor, we show that there exists one dominating $B$. When this divisor is potentially dense and the geometric generic fibre of the Lagrangian fibration is simple, we can use the divisor to construct a nondegenerate section of the pulled-back abelian scheme, which implies potential density of rational points on $X$.

Finally, to prove the main theorem, we construct a projective Shafarevich--Tate twist $Y$ of $X$ admitting a dominant rational map to $X$. We then show that $Y$ admits two distinct rational Lagrangian fibrations. One fibration is inherited from the original fibration on $X$. When $Y$ has Picard rank $2$, we choose the twist so that it has no algebraic MBM classes by \cite{Amerik_Verbitsky_20}. Then by \cite{Amerik_Verbitsky_15}, the other isotropic boundary ray is nef. When $Y$ has higher Picard rank, we use lattice theory to construct a rational isotropic vector whose reflections in prime exceptional divisors give a second movable isotropic ray by \cite{Riess_23}. In both cases, the second ray induces a rational Lagrangian fibration, and the two-fibration theorem implies potential density of rational points on $Y$, hence of $X$.

\subsection{Notations and terminology}

Throughout this manuscript, we write $\rho(X)=\rank\NS(X)$ for the Picard rank of $X$, and $\rho(X_{\bar{k}})$ for its geometric Picard rank. We write $b_n(X)=\dim_{\QQ}H^n(X(\CC),\QQ)$ for its $n$th Betti number. Furthermore, we denote by $q_X$ the Beauville--Bogomolov quadratic form on $H^2(X,\ZZ)$, and by $(\cdot,\cdot)$ its associated bilinear form, normalized so that $q_X(x+y)=q_X(x)+2(x,y)+q_X(y)$.

For an integral lattice $\Lambda$ and a nonzero vector $v\in\Lambda$, we denote by $\operatorname{div}_{\Lambda}(v)$ the positive generator of the ideal $(v,\Lambda)\subset\ZZ$. When the ambient lattice is clear, we simply write $\operatorname{div}(v)$. Furthermore, the notation $\pr_A(v)$ denotes the orthogonal projection of a vector $v$ onto a nondegenerate subspace $A$ with respect to the given bilinear form. 

\section*{Acknowledgement and AI statement}
The author would like to thank ICERM for its hospitality during the conference ``Cubic Fourfolds, Gushel--Mukai Fourfolds, and Hyperk\"ahler Manifolds.'' The author would also like to thank Brendan Hassett for helpful discussions, Zhizhong Huang for suggestions on the problem and terminology, and Bianca Viray and Max Lieblich for their help and support.

The question studied in this paper was initially raised during a talk by Zhizhong Huang at a conference in December 2025, concerning whether hyperk\"ahler varieties admitting a Lagrangian fibration satisfy potential Hilbertian properties. Subsequently, Zhizhong Huang suggested that the author consider potential density of the rational points on these varieties, which might be more approachable and require fewer dynamical arguments. The author later came across the work of Hassett--Tschinkel \cite{Hassett_Tschinkel_00} and observed that their approach could be used to establish the potential density property in the case of two distinct Lagrangian fibrations.

Subsequently, the author explored this problem using AI and received an output of approximately 500 pages, containing purported proofs for $K3^{[4]}$-type varieties, generalized Kummer fourfolds, OG6, and OG10, respectively. It also contained an unverified proof of Corollary~\ref{thm:witt-two-density} using moduli-theoretic methods. After examining some of these arguments, the author realized that the results could be proved uniformly using the approach developed in Section~4. The remaining mathematical value of the 500-page document is unclear. It may contain further interesting results, but a complete examination is beyond the author's current capacity. In particular, Sections~3 and~4 contain arguments partially suggested by AI, while the Appendix provides an example of a typical argument generated by AI.

Overall, the paper was written entirely by the author and subsequently checked and reviewed with the assistance of AI. The author takes full responsibility for the contents of the paper, including the correctness of all mathematical statements and proofs.

\section{Hyperkählers admitting two Lagrangian fibrations}
In this section, after establishing some preliminary results, we show that rational points on irreducible hyperk\"ahler varieties over number fields admitting two distinct rational Lagrangian fibrations are potentially dense.

\subsection{Preliminary results}
In this subsection, we state and prove some basic technical results which will be used later in our proofs. We first use a multisection to construct a finite \'etale morphism from a torsor under an abelian scheme to the abelian scheme, whose degree is determined by that of the multisection. We then prove a density result for a dominant morphism whose generic fibre is an abelian torsor, using an auxiliary variety admitting a rational map to its total space. Next, we show that a section is nondegenerate at the geometric generic point if it is nondegenerate at a closed point of the base. Finally, we show that every positive-dimensional subvariety of an abelian variety contains a nontorsion algebraic point in each nonempty open subset.

\begin{lemma}\label{lem:auxiliary-multiplication}
    Let $f: T\to B$ be a torsor under an abelian scheme $A\to B$, where $B$ is geometrically integral over a field of characteristic $0$. After replacing $B$ by a nonempty open subset, there is a finite \'etale multisection $Z\subset T$ of degree $d>0$ and a finite \'etale surjective $B$-morphism $\varphi: T\rightarrow A$ with $\varphi(t) = \sum_{z\in Z_{f(t)}} (t-z)$. If the relative dimension of $T/B$ is $n$, then $\deg\varphi=d^{2n}$.
\end{lemma}

\begin{proof}
    Let $\eta$ be the generic point of $B$, and let $p$ be a closed point of $T_\eta$. Assume that the residue field $k(p)$ has degree $d$ over $k(B)$. Then, after shrinking $B$, the closure of $p$ in $T$ gives a finite \'etale multisection $Z\to B$ of degree $d$. The map $T\times_B Z\to A$ given by $(t,z)\mapsto t-z$ is a $B$-morphism. Applying the trace map along the finite \'etale morphism $T\times_B Z\to T$ gives the desired morphism $\varphi:T\to A$. By construction, for every $u\in A_{f(t)}$, we have $\varphi(t+u)=\varphi(t)+[d]u$. Hence, after passing to an \'etale cover of $B$ and choosing an origin, $\varphi$ is multiplication by $d$ followed by a translation. Therefore, $\varphi$ is finite \'etale and surjective of degree $d^{2n}$.
\end{proof}

\begin{proposition}\label{prop:parameter-core}
    Let $f:Y\to B$ be a dominant morphism whose generic fibre is a torsor under the generic fibre of an abelian scheme $A\to B$. Let $M$ be a geometrically integral variety admitting a rational map $j:M\dashrightarrow Y$ dominating $B$, such that $M(k)$ is Zariski dense. Let $s$ be a rational section of $A_M\to M$ which is nondegenerate at the geometric generic point of $M$. Then $Y(k)$ is Zariski dense.
\end{proposition}

\begin{proof}
    The argument is the same as that of \cite[Proposition~3.4]{Hassett_Tschinkel_00}. After shrinking $M$ and $B$, we may assume that $f:Y\to B$ is a torsor under $A\to B$ and that $j$, $h=f\circ j$, and $s$ are morphisms. Consider the torsor $Y_M=Y\times_B M\to M$ under $A_M$. Then $j$ induces a section $e:M\to Y_M$. For each $r\in\ZZ$, let $e_r=e+[r]s$. Since $M(k)$ is Zariski dense in $M$, the rational points $e_r(M(k))$ are Zariski dense in $e_r(M)$.

    Let $\bar\eta$ be a geometric generic point of $M$. By the nondegeneracy of $s(\bar\eta)$, the points $e(\bar\eta)+[r]s(\bar\eta)$, for $r\in\ZZ$, are Zariski dense in the geometric generic fibre of $Y_M\to M$. Hence the union of all $e_r(M)$ is Zariski dense in $Y_M$. Therefore, $Y_M(k)$ is Zariski dense in $Y_M$. Since $h:M\to B$ is dominant, the projection $Y_M\to Y$ is dominant. Hence $Y(k)$ is Zariski dense in $Y$. 
\end{proof}

\begin{lemma}\label{lem:cyclic-specialization}
    Let $f:A\to S$ be an abelian scheme over a smooth geometrically integral variety $S$ over a field of characteristic $0$, and let $s:S\to A$ be a section. If $s(x)$ is nondegenerate at some closed point $x\in S$, then $s$ is nondegenerate at the geometric generic point of $S$.
\end{lemma}

\begin{proof}
    After a finite extension and shrinking $S$ around $x$, we may assume that $x$ is a rational point and $S$ is quasi-projective. Consider the collection of sections $E=\{[m]s:m\in\ZZ\}$. Since $s(x)$ is nondegenerate, the set $E(x)=\{[m]s(x):m\in\ZZ\}$ is Zariski dense in $A_x$. By \cite[Lemma~4.1]{Bogomolov_Halle_Halvard_Tanimoto_18}, the union $\bigcup_{m\in\ZZ}[m]s(S)$ is Zariski dense in $A$.

    Let $\eta$ be the generic point of $S$. Suppose that $s$ is degenerate at the geometric generic point. Since all its multiples are defined over $k(S)$, the geometric Zariski closure of the cyclic subgroup descends to a proper closed subgroup $H_\eta\subset A_\eta$. Let $H$ be its schematic closure in $A$. Then $H$ is a proper closed subset of $A$ containing the images of all sections $[m]s$, contradicting the Zariski density of $\bigcup_{m\in\ZZ}[m]s(S)$. Hence $s$ is nondegenerate at the geometric generic point.
\end{proof}

\begin{lemma}\label{lem:abelian_subvar_nontorsion}
    Let $A$ be an abelian variety over a number field $k$. Then every positive-dimensional closed integral subvariety $Z\subset A$ defined over $k$ contains a nontorsion algebraic point in each nonempty open subset of $Z$.
\end{lemma}

\begin{proof}
    Let $W\subset Z$ be a nonempty open subset. After a finite extension of $k$, choose an integral projective curve $C\subset Z$ such that $C\cap W\neq\varnothing$. Let $\widetilde{C}$ be its normalization, and let $\nu:\widetilde{C}\to C$ be the normalization map. Then choose a nonconstant finite morphism $\pi:\widetilde{C}\to\PP^1$ of degree $d>0$. Taking the fibres above infinitely many points of $\PP^1(k)$, and excluding finitely many fibres meeting $\widetilde{C}\setminus\nu^{-1}(W)$, we obtain infinitely many distinct algebraic points of $W$ of degree at most $d$ over $k$.

    Let $\cL$ be a symmetric ample line bundle on $A$. Then the canonical height $\widetilde{h}_{\cL}$ vanishes on torsion points and differs from the usual Weil height $h_{\cL}$ by a bounded function \cite[\S B.5]{Hindry_Silverman_00}. By Northcott's theorem \cite[Theorem~B.2.3]{Hindry_Silverman_00}, there are only finitely many torsion points of bounded degree over $k$. Therefore, at least one of the algebraic points constructed above is nontorsion.
\end{proof}

\subsection{Rational Lagrangian fibrations and potential density for two fibrations}
There are different usages of the term \emph{rational Lagrangian fibration} in the literature, as noted in \cite{Kim_Laza_Martin_26}. We follow the terminology of \cite{Markushevich_06, Kamenova_Lehn_25}.

\begin{definition}
    A \emph{rational Lagrangian fibration} of $X$ is a rational map $f\colon X\dashrightarrow B$ such that there exist a smooth projective hyperk\"ahler variety $X'$, a birational map $h\colon X\dashrightarrow X'$, and a regular Lagrangian fibration $g\colon X'\to B$ satisfying $f=g\circ h$.
\end{definition}

We now present a lemma showing that a rational Lagrangian fibration can be realized over a finite extension of a number field.

\begin{lemma}\label{lem:descent_finite_extension}
    Let $X/k$ be an irreducible hyperk\"ahler variety over a number field $k$. Assume that $X_\CC$ admits a rational Lagrangian fibration. Then $X$ admits a rational Lagrangian fibration after a finite extension of $k$.
\end{lemma}

\begin{proof}
    Let $Y_\CC$ be a hyperk\"ahler variety birational to $X_\CC$, admitting a regular Lagrangian fibration $f:Y_\CC\to B_\CC$. Choose a symplectic form $\sigma$ on $Y_\CC$ and a smooth point $b\in B_\CC$ such that $f^{-1}(b)$ is smooth and connected. Choose a finitely generated extension $K/\bar{k}$ over which all these data, together with the birational maps between $Y_\CC$ and $X_\CC$ and their inverses, are defined. We can then spread them out over a smooth integral scheme $S$ of finite type over $\bar{k}$ with function field $K$.

    After shrinking $S$, we may assume that $\mathcal Y\to S$ is smooth and projective with geometrically integral fibres, and $\mathcal B\to S$ is flat and projective with geometrically integral and normal fibres of dimension $n$. The birational maps between $\mathcal Y$ and $X_{\bar{k}}\times S$, together with their inverses, are defined on dense open subsets of every fibre, and their compositions remain the identity on these opens. Moreover, $\sigma$ extends to a relative $2$-form on $\mathcal Y$ which is nondegenerate, and $f$ extends to a morphism $F:\mathcal Y\to\mathcal B$. The point $b$ extends to a section of the smooth locus of $\mathcal B\to S$, and its corresponding family of fibres remains smooth and geometrically connected. We may also assume that the relative symplectic form vanishes on the tangent spaces of the smooth fibres of $F$.

    We now show that the fibres of each specialized morphism $F_s:\mathcal Y_s\to\mathcal B_s$ are geometrically connected. Since the fibre over $b_s$ is smooth and geometrically connected, there exists a nonempty open subset $V_s\subset\mathcal B_s$ containing $b_s$ such that $F_s^{-1}(V_s)\to V_s$ is smooth and proper. By \cite[Tag~0E0N]{stacks-project}, the number of geometric connected components of the fibres of a smooth proper morphism is locally constant. Since $V_s$ is integral and contains $b_s$, the geometric generic fibre of $F_s$ is connected. Then by Stein factorization and normality of $\mathcal B_s$ \cite[Tag~0AY8]{stacks-project}, all geometric fibres of $F_s$ are connected.

    Choose a closed point $s\in S(\bar{k})$. Then $\mathcal Y_s$ is birational to $X_{\bar{k}}$ and is hyperk\"ahler, since its complexification is deformation equivalent to $Y_\CC$ in a smooth projective family. Furthermore, $F_s:\mathcal Y_s\to\mathcal B_s$ is surjective, since it is smooth over the nonempty open subset $V_s$ and proper. Its geometric generic fibre has dimension $n$ and is Lagrangian with respect to the specialized symplectic form. Thus $F_s$ is a regular Lagrangian fibration. Since all these data are of finite presentation over $\bar{k}$, they descend to a finite extension of $k$.
\end{proof}

Let $f:X\dashrightarrow B$ be a rational Lagrangian fibration, with $g$, $h$, and $X'$ as in the preceding definition. Let $H$ be an ample divisor on $B$. Note that $\rho(B_{\bar{k}})=1$ by \cite[Theorem~2.10]{Kamenova_Lehn_25}. Let $D=h^*g^*H\in\operatorname{NS}(X_{\bar{k}})_{\mathbb Q}.$ Then $D$ is a nonzero movable isotropic divisor class; see, for example, \cite[Remark~2.12]{Kamenova_Lehn_25}.

\begin{definition}
    We call the ray $\mathbb R_{\geq 0}[D] \subset\operatorname{NS}(X_{\bar{k}})_{\mathbb R}$ the \emph{movable isotropic divisor ray associated with} the rational Lagrangian fibration $f:X\dashrightarrow B$.
\end{definition}

We first describe the fibres of the joint rational map associated with two Lagrangian fibrations. Using this description, we then show that two distinct movable isotropic divisor rays give a dominant generically finite rational map to the product of their bases.

\begin{lemma}\label{lem:joint-fibers}
    Let $X$ be an irreducible hyperk\"ahler variety of dimension $2n$ over a field of characteristic zero. Let $f:X\to B_1$ be a regular Lagrangian fibration and $g:X\dashrightarrow B_2$ be a rational Lagrangian fibration. Let $Y$ be the closure of the image of $(f,g)$. Suppose that $\dim Y=n+r<2n$. Then the closures of the irreducible components of the geometric general fibres of $(f,g)$ are translates of abelian subvarieties of dimension $n-r$ in the smooth fibres of $f$.
\end{lemma}

\begin{proof}
    We may work over an algebraically closed field, since the assertion concerns geometric fibres. Let $(b,y)$ be a general geometric point of $Y$, and let $A=f^{-1}(b)$ be the smooth fibre of $f$. Take smooth open subsets of the two bases with \'etale coordinates $t_1,\ldots,t_n$ and $s_1,\ldots,s_n$, such that on a common open subset $U$ of $X$, both fibrations are defined and smooth. Then at $(b,y)$, the conormal space of $Y$ in $B_1\times B_2$ has dimension $n-r$. Let $G_1,\ldots,G_{n-r}$ be local equations defining $Y$ near this point and write
    $$dG_\ell=\sum_{i=1}^n a_{\ell i}\,dt_i+\sum_{j=1}^n b_{\ell j}\,ds_j.$$
    Since the projections $Y\to B_i$ are dominant and generically smooth, the projections of the conormal space to the two cotangent spaces are injective. Hence the vectors $(a_{\ell 1},\ldots,a_{\ell n})$, for $1\leq\ell\leq n-r$, are linearly independent.
    
    Let $\sigma$ be the symplectic form on $X$, and let $X_u$ be the Hamiltonian vector field associated with a local function $u$, i.e. $\iota_{X_u}\sigma=du$. Then the Lagrangian property implies that $X_{f^*t_1},\ldots,X_{f^*t_n}$ are tangent to the fibres of $f$ and form a frame of their tangent bundles. Similarly for $X_{g^*s_1},\ldots,X_{g^*s_n}$. Differentiating $G_\ell(f,g)=0$ and applying $\sigma^{-1}$, we have
    $$V_\ell:=\sum_{i=1}^n a_{\ell i}(f,g)X_{f^*t_i} =-\sum_{j=1}^n b_{\ell j}(f,g)X_{g^*s_j}.$$
    Hence the $V_\ell$ are tangent to both fibrations. Moreover, by the linear independence of $(a_{\ell 1},\ldots,a_{\ell n})$ for $1\leq\ell\leq n-r$, the vector fields $V_\ell$ are linearly independent.

    After choosing an origin, $A$ is an abelian variety, and the vector fields $X_{f^*t_i}|_A$ are translation invariant. Let $E$ be the reduced closure in $A$ of an irreducible component of the fibre of $(f,g)$ over $(b,y)$ on $U$. Then $\dim E=n-r$. Along $E\cap U$, $(f,g)$ is constant. Thus the $a_{\ell i}(f,g)$ are constant on $E\cap U$, and hence the $V_\ell$ give $n-r$ independent translation-invariant vector fields on $A$ tangent to $E$ on a dense open subset. Since $E$ is reduced, these vector fields are tangent to all of $E$. Therefore, their span is contained in the Lie algebra of the translation stabilizer $\operatorname{Stab}_A(E)$. Hence
    $$\dim\operatorname{Stab}_A(E)^\circ\geq n-r=\dim E.$$
    Since the action of $A$ on itself is free, we also have $\dim\operatorname{Stab}_A(E)^\circ\leq\dim E$. Hence equality holds, and $E$ is a translate of $\operatorname{Stab}_A(E)^\circ$.
\end{proof}

\begin{proposition}\label{prop:transversality}
    Let $X$ be an irreducible hyperk\"ahler variety of dimension $2n$ over a field of characteristic zero. Let $f:X\to B_1$ be a regular Lagrangian fibration and $g:X\dashrightarrow B_2$ be a rational Lagrangian fibration. Denote their associated movable isotropic divisor rays by $R_1$ and $R_2$.

    If $R_1\neq R_2$, then the rational map $(f,g):X\dashrightarrow B_1\times B_2$ is dominant and generically finite.
\end{proposition}

\begin{proof}
    We may work over $\CC$, after passing to a finitely generated field of definition if necessary, since dominance and generic finiteness can be checked after base change. Let $D_1,D_2$ be nonzero divisor classes in the rays $R_1$ and $R_2$, respectively, chosen as pullbacks of ample divisors from the two bases on their regular models. Since these rays are distinct and lie in the closure of the same positive cone, the Beauville--Bogomolov form gives $q_X(D_1,D_2)>0$ and $q_X(D_1)=q_X(D_2)=0$ by \cite[Lemma~2.5]{Kamenova_Lehn_25}. Then, applying Fujiki's relation to $sD_1+tD_2$ and comparing the coefficients, we have
    \begin{equation}\label{eq:transversality-intersection}
    \binom{2n}{n}\int_X D_1^nD_2^n=c_X2^n q_X(D_1,D_2)^n>0.
    \end{equation}

    Let $A$ be a smooth general fibre of $f$. We show that $g|_A$ is nonconstant. Suppose the contrary. Since the birational map $h:X\dashrightarrow X'$, where $X'$ is the regular model of $g$, is an isomorphism in codimension one by \cite[Remark~3.4]{Looijenga_21}, the map $g$ is defined outside a subset of codimension at least two in $A$. Since $g|_A$ is constant, $D_2|_A$ is the pullback of a divisor class from a point on this open subset, hence trivial. Since $A$ is smooth, $D_2|_A$ is trivial on all of $A$ by Hartogs' theorem. As $D_1=f^*H_1$ for an ample divisor $H_1$ on $B_1$, the projection formula gives
    $$\int_X D_1^nD_2^n =\left(\int_{B_1}H_1^n\right) \left(\int_A (D_2|_A)^n\right)=0,$$
    contradicting \eqref{eq:transversality-intersection}. Therefore, $g|_A$ is nonconstant.

    Let $Y$ be the closure of the image of $(f,g)$. Since the projections $Y\to B_i$ are dominant and $g|_A$ is nonconstant, we have $\dim Y=n+r$ for some $1\leq r\leq n$. We show that $\dim Y=2n$ by contradiction. Assume that $0<r<n$. Then by Lemma~\ref{lem:joint-fibers}, the closures of the irreducible components of the geometric generic fibre of $(f,g)$ are translates of abelian subvarieties of dimension $n-r$. In particular, their Kodaira dimension is $0$.

    After taking a resolution and Stein factorization of the joint rational map $(f,g):X\dashrightarrow Y$, its geometric generic fibre is birational to one of these components, and hence has Kodaira dimension $0$. Then by \cite[Proposition~3.4]{Amerik_Campana_08}, its dimension is at least $n$, contradicting $n-r<n$. Therefore, $\dim Y=2n=\dim X=\dim(B_1\times B_2)$. Hence $(f,g)$ is dominant and generically finite.
\end{proof}

\begin{corollary}\label{cor:dominant-two-rational}
    Let $k$ be a number field and $X/k$ be an irreducible hyperk\"ahler variety of dimension $2n$. Assume that $X$ admits two rational Lagrangian fibrations $f_i:X\dashrightarrow B_i$ whose associated movable isotropic divisor classes span distinct rays. Then $(f_1,f_2):X\dashrightarrow B_1\times B_2$ is dominant and generically finite.
\end{corollary}

\begin{proof}
    After replacing $k$ with a finite extension $K$, we may assume that the regular models of $f_{i,K}:X_K\dashrightarrow B_{i,K}$ are defined over $K$. We omit the base change notation in the rest of the proof. Let $X_1,X_2$ be the hyperk\"ahler models for the two maps, respectively, with birational maps $h_i\colon X\dashrightarrow X_i$ and regular Lagrangian fibrations $\varphi_i\colon X_i\to B_i$, such that $f_i=\varphi_i\circ h_i$. Consider the map $\psi:=\varphi_2\circ h_2\circ h_1^{-1}\colon X_1\dashrightarrow B_2.$     Then $\psi$ is a rational Lagrangian fibration. Furthermore, since birational maps between hyperk\"ahler varieties induce isometries of their Beauville--Bogomolov lattices and identify their movable cones, the corresponding movable isotropic divisor rays on $X_1$ are still distinct. Hence $\varphi_1$ and $\psi$ satisfy the assumptions of Proposition~\ref{prop:transversality}. Therefore, $(\varphi_1,\psi)$ is dominant and generically finite, and so is $(f_1,f_2)=(\varphi_1,\psi)\circ h_1.$ Since dominance and generic finiteness can be checked after base change, the result follows over $k$.
\end{proof}

Now we can prove the theorem of this section. When $X$ admits two rational Lagrangian fibrations whose associated movable isotropic divisor rays are distinct, Corollary~\ref{cor:dominant-two-rational} supplies a dominant and generically finite map to the product of their bases. This allows us to construct a potentially dense multisection of one fibration from a smooth fibre of the other. The preliminary results then give a nondegenerate section of the pulled-back abelian scheme, allowing us to deduce potential density of rational points on $X$.

\begin{theorem}\label{thm:two-fibrations}
    Let $k$ be a number field and $X/k$ be an irreducible hyperk\"ahler variety of dimension $2n$. Assume that $X$ admits two rational Lagrangian fibrations $f_i:X\dashrightarrow B_i$ whose associated movable isotropic divisor classes span distinct rays. Then the rational points on $X$ is potentially dense.
\end{theorem}

\begin{proof}
    After replacing $k$ with a finite extension $K$, we may assume that the regular models of $f_{i,K}:X_K\dashrightarrow B_{i,K}$ are defined over $K$. Write $f_i=\varphi_i\circ h_i$, where $h_i:X\dashrightarrow X_i$ are birational maps and $\varphi_i:X_i\to B_i$ are regular Lagrangian fibrations. We omit the base change notation in the rest of the proof. By Corollary~\ref{cor:dominant-two-rational}, $(f_1,f_2)$ is dominant and generically finite.

    By shrinking $B_i$ to nonempty open subsets $B_i'$, we may assume that $\varphi_i^{-1}(B_i')\to B_i'$ is a torsor $T_i\to B_i'$ under an abelian scheme $A_i\to B_i'$, see \cite[Lemma~1.16]{Huybrechts_Mauri_22} and \cite[Theorem~5.3]{Laurent_Schroer_24}. Choose a common open subset $U$ of $X,X_1,X_2$ such that both fibrations are regular on $U$ and their images lie in $B_i'$. Then $F=(f_1,f_2)|_U\colon U\to B_1'\times B_2'$ is dominant. By Chevalley's theorem \cite[Exercise~II.3.19]{Hartshorne_77}, $F(U)$ is constructible and hence contains a nonempty open subset $W\subset B_1'\times B_2'$. Replacing $U$ by $F^{-1}(W)$, we obtain a surjective morphism $F\colon U\to W$.

    By Lemma~\ref{lem:auxiliary-multiplication}, after shrinking $B_1'$ and $W$ accordingly, there is a finite \'etale surjective morphism $\varphi:T_1\to A_1$ whose restriction to every geometric fibre is multiplication by $d$ followed by translation, after choosing an origin. Choose a rational point $b$ in the image of $W\to B_1'$ after a finite extension. Then $U_b$ is a nonempty open subset of $T_{1,b}$. Since $\varphi_b:T_{1,b}\to A_{1,b}$ is finite \'etale and surjective, $\varphi_b(U_b)\subset A_{1,b}$ is a nonempty open subset. By \cite[Proposition~3.1]{Hassett_Tschinkel_00}, after a finite extension, there exists $a\in A_{1,b}$ whose cyclic subgroup is geometrically Zariski dense. Hence we can choose a nonzero multiple $q=[m]a\in\varphi_b(U_b)$, which is still nondegenerate. After another finite extension, we can choose $p\in U_b$ such that $\varphi_b(p)=q$.

    Let $\pi_i:U\to B_i'$ be the restrictions of the two fibrations. Let $y=\pi_2(p)$ and $M=\varphi_2^{-1}(y)$. Then $M$ is a torsor under an abelian variety and contains $p$. Choosing $p$ as the origin, $M$ is an abelian variety. We show that $M$ dominates $B_1'$. Since $(b,y)\in W$, the fibre $W_y$ is a nonempty open subset of $B_1'$. Since $F$ is surjective, for every $x\in W_y$, there exists $w\in U$ such that $F(w)=(x,y)$. Hence $w\in M\cap U$, and $x$ lies in the image of $\pi_1|_{M\cap U}:M\cap U\to B_1'$. Therefore, this map is dominant. Since $\dim M=\dim B_1'=n$, it is also generically finite. Hence $M\cap U$ defines a multisection of the first fibration.

    Since $M$ is an abelian variety, by \cite[Proposition~3.1]{Hassett_Tschinkel_00}, it has a nondegenerate rational point after a finite extension. Hence rational points on $M\cap U$ is potentially dense. The inclusion $j:M\cap U\to T_1$ induces a section $e:M\cap U\to T_{1,M\cap U}$. Composing $e$ with the pullback of $\varphi$ gives a section $s:M\cap U\to A_{1,M\cap U}$, where $s(p)=\varphi_b(p)=q$. Since $q$ is nondegenerate, Lemma~\ref{lem:cyclic-specialization} implies that $s$ is nondegenerate at the geometric generic point of $M\cap U$.

    Applying Proposition~\ref{prop:parameter-core} to $T_1\to B_1'$ with $M\cap U$ as the parameter space, we see that $T_1$, and hence rational points on $X_1$, is potentially dense. Since potential density is a birational invariant, rational points on $X$ is also potentially dense.
\end{proof}

\section{Hyperkählers admitting one Lagrangian fibration}
In this section, we record some general results related to Conjecture~\ref{conj:one-fibration}.

\subsection{Potential density with no prime exceptional divisors}
In this subsection, we prove a general potential density result for rational points on hyperk\"ahler varieties admitting a Lagrangian fibration and having no prime exceptional divisors, under conditions depending on their Picard rank.

Note that the existence of a Lagrangian fibration implies $\rho(X_{\bar{k}})\geq 2$, since its associated isotropic divisor class cannot be proportional to an ample class. We first review the terminology introduced in \cite{Markman_13}.

\begin{definition}
    [\cite{Markman_13}] A divisor $E$ on $X$ is called a \emph{prime exceptional divisor} if $E$ is reduced, irreducible, and has negative Beauville--Bogomolov square, i.e. $q_X(E)<0$. 

    For every prime exceptional divisor $E$ on $X$, we define its \emph{reflection} by
    $$R_E(\alpha)=\alpha-\frac{2(\alpha,[E])}{q_X([E])}[E],$$
    and its \emph{exceptional reflection group} by
    $$W_{\operatorname{Exc}}(X)= \left\langle R_E : E\subset X \text{ is a prime exceptional divisor}\right\rangle.$$
\end{definition}

Using Theorem~\ref{thm:two-fibrations}, we can quickly prove the case $\rho(X_{\bar{k}})\geq 3$.

\begin{theorem}\label{thm:fibration_rank_3}
    Let $X$ be an irreducible hyperk\"ahler variety over a number field $k$, admitting a regular Lagrangian fibration and having no geometric prime exceptional divisors. Assume that $\rho(X_{\bar{k}})\geq 3$. Then rational points on $X$ is potentially dense.
\end{theorem}

\begin{proof}
    First assume that $\rho(X_{\bar{k}})\geq 3$, and let $R$ be the movable isotropic divisor ray associated with the fibration $f:X\to B$. Since $X_{\bar{k}}$ has no prime exceptional divisors, its exceptional reflection group $W_{\operatorname{Exc}}(X_{\CC})$ is trivial. Then by \cite[Theorem~6.18 and Lemma~6.23]{Markman_11}, the image of the birational automorphism group $\operatorname{Bir}(X_{\CC})$ in $O^+(\NS(X_{\CC}))$ has finite index. By the Borel density theorem \cite{Borel_60}, its intersection with $\mathrm{SO}(\NS(X_{\CC})_\RR)$ is Zariski dense in $\mathrm{SO}(\NS(X_{\CC})_\RR)$. Since the stabilizer of $R$ is a proper algebraic subgroup of $\mathrm{SO}(\NS(X_{\CC})_\RR)$, there exists $\varphi\in\operatorname{Bir}(X_{\CC})$ such that $\varphi^*R \neq R$. After replacing $k$ with a finite extension, we may assume that $\varphi$ is defined over $k$. Then $f\circ\varphi$ is a second rational Lagrangian fibration with associated ray $\varphi^*R$. By Theorem~\ref{thm:two-fibrations}, rational points on $X$ is potentially dense.
\end{proof}

\begin{remark}
    Note that the Picard rank condition is important here, since the identity component of $\mathrm{SO}(\NS(X_{\CC})_\RR)$ is noncompact and semisimple without compact factors when $\rho(X_{\bar{k}})\geq 3$, as the Beauville--Bogomolov form has signature $(1,\rho-1)$. If the Picard rank is $2$, then the identity component of $\mathrm{SO}(\NS(X_{\CC})_\RR)$ is isomorphic to $\RR_{>0}$, which is not semisimple and does not give enough room to obtain another rational Lagrangian fibration by the same argument.
\end{remark}

We now consider the case $\rho(X_{\bar{k}})=2$ with $b_2(X)\geq 5$. We will again obtain a second rational Lagrangian fibration, this time by showing that the other isotropic boundary ray of the positive cone becomes nef on a smooth birational model. The semiampleness theorem of Engel--Mauri \cite{Engel_Mauri_26} then gives the required fibration.

\begin{lemma}\label{lem:finite_negative_lattice}
    Let $\Lambda$ be an integral lattice of signature $(1,1)$ containing a nonzero rational isotropic vector. Then for every $M>0$, the set $\{z\in\Lambda\mid -M\leq q(z)<0\}$ is finite. 
\end{lemma}

\begin{proof}
    Since $\Lambda$ contains a nonzero rational isotropic vector, it is a split quadratic plane over $\QQ$. Let $l,r\in\Lambda$ be generators of its two isotropic lines, and let $c=q(l,r)\neq 0$. Since $l,r$ span $\Lambda_\QQ$, there is a positive integer $d$ such that for any $z\in\Lambda$, we have $z=\frac{a}{d}l+\frac{b}{d}r$, where $a,b\in\ZZ$. Since $l,r$ are isotropic, we have $q(z)=\frac{2cab}{d^2}$. If $-M\leq q(z)<0$, then $0<|ab|\leq\frac{Md^2}{2|c|}.$ Hence there are only finitely many pairs of nonzero integers $(a,b)$, and therefore only finitely many such $z\in\Lambda$.
\end{proof}

We first recall some definitions from \cite{Amerik_Verbitsky_15,Amerik_Verbitsky_20}. The condition $b_2(X)\geq 5$ is required by the boundedness result for MBM classes that we use later.

\begin{definition}
    A nonzero class $z\in\NS(X)$ with $q_X(z)<0$ is called \emph{monodromy birationally minimal (MBM)} if there exist a Hodge monodromy operator $\gamma$ and a smooth hyperk\"ahler variety $X'$ birational to $X$ such that $(\gamma z)^\perp$ contains a codimension-one face of the pullback of the K\"ahler cone of $X'$.

    A \emph{K\"ahler--Weyl chamber} of $X$ is the image of the K\"ahler cone of some birational model $X'$ under a Hodge monodromy operator $\gamma$, after identifying $H^2(X',\RR)$ with $H^2(X,\RR)$ via the birational map.
\end{definition}

\begin{proposition}\label{prop:other_ray_nef}
    Let $X$ be an irreducible hyperk\"ahler variety over $\CC$ with $\rho(X)=2$ and $b_2(X)\geq 5$. Assume that $\NS(X)$ contains a nonzero isotropic class $l$ in the closure of the positive cone, and that $X$ has no prime exceptional divisors. Then the other isotropic boundary ray of the positive cone in $\NS(X)_\RR$ is rational and becomes nef on a smooth birational hyperk\"ahler model of $X$.
\end{proposition}

\begin{proof}
    Let $h$ be an ample class in $\NS(X)$. Let $r=h-\frac{q_X(h)}{2(l,h)}l$. Then by direct computation, we have 
    $$q_X(r)=0,\qquad (r,h)=\frac{q_X(h)}{2}>0,\qquad (r,l)=(h,l)>0.$$
    Hence $r$ generates the other rational isotropic boundary ray of the positive cone.

    By \cite[Theorem~3.17]{Amerik_Verbitsky_20}, the negative squares of primitive MBM classes are bounded when $b_2(X)\geq 5$. Then since $\NS(X)$ has signature $(1,1)$ and contains a nonzero rational isotropic class, Lemma~\ref{lem:finite_negative_lattice} shows that there are finitely many primitive algebraic MBM classes, hence finitely many MBM hyperplanes in the positive cone. Furthermore, by the signature of $\NS(X)$, the positive cone is the sector between the two isotropic rays generated by $l$ and $r$. For any class $z$ of negative square, its orthogonal hyperplane cannot contain $r$, since $r^\perp=\RR r$ contains no vectors of negative square. Hence each MBM hyperplane intersects the positive cone in an interior ray distinct from $r$. Since there are only finitely many such hyperplanes, one of the chambers of their complement is adjacent to the boundary ray $\RR_{\geq 0}r$. Hence its closure contains $r$.
    
    By \cite[Theorem~6.2]{Amerik_Verbitsky_15}, these chambers are K\"ahler--Weyl chambers. Since $X$ has no prime exceptional divisors, $W_{\operatorname{Exc}}(X)$ is trivial by \cite[Theorem~6.18(3)]{Markman_11}. Hence the Hodge monodromy group coincides with the birational monodromy group by \cite[Theorem~6.18(5)]{Markman_11}. Therefore, the chamber whose closure contains $r$ is the pullback of the K\"ahler cone of a smooth hyperk\"ahler variety $Y$ birational to $X$. Hence $r$ is nef on $Y$. Furthermore, the chamber contains a rational class of positive square, whose positive integral multiple is ample on $Y$. Hence $Y$ is projective.
\end{proof}

Combining Proposition ~\ref{prop:other_ray_nef} with the semiampleness theorem of Engel--Mauri \cite{Engel_Mauri_26}, we now obtain a second rational Lagrangian fibration and prove potential density in Picard rank two.

\begin{theorem}\label{thm:fibration_rank_2}
    Let $X$ be an irreducible hyperk\"ahler variety over a number field $k$, admitting a regular Lagrangian fibration and having no geometric prime exceptional divisors. Assume that $\rho(X_{\bar{k}})=2$ and $b_2(X)\geq 5$. Then rational points on $X$ is potentially dense.
\end{theorem}

\begin{proof}
    Assume that $\rho(X_{\bar{k}})=2$. By Proposition~\ref{prop:other_ray_nef}, the other isotropic boundary ray becomes nef on a smooth hyperk\"ahler birational model $X'$ of $X_{\CC}$. Since this ray is rational, it contains the class of a nontrivial nef isotropic line bundle. By the hyperk\"ahler SYZ theorem \cite[Theorem~1.1]{Engel_Mauri_26}, this line bundle is semiample, hence defines a regular Lagrangian fibration on $X'$ whose associated ray is distinct from $R$. After replacing $k$ with a finite extension, we may assume that this second rational Lagrangian fibration is defined over $k$. Applying Theorem~\ref{thm:two-fibrations}, we obtain potential density of rational points on $X$.
\end{proof}

\subsection{Potential density from that of prime exceptional divisors}

We now show that potential density of rational points on prime exceptional divisors implies potential density of rational points on the ambient variety, under the same Picard rank conditions. Together with the preceding results for varieties without prime exceptional divisors, this gives a general criterion for potential density.

\begin{lemma}\label{lem:maximal-rank}
    Let $X$ be an irreducible hyperk\"ahler variety over $\CC$ with $\rho(X)\geq 3$. If $X$ has a prime exceptional divisor, then the classes of its prime exceptional divisors span $\NS(X)_\QQ$.
\end{lemma}

\begin{proof}
    Let $R\subset\NS(X)_\QQ$ be the subspace generated by the classes of prime exceptional divisors. For every prime exceptional divisor $D$, consider its reflection $R_D(x)=x-\frac{2(x,[D])}{q_X(D)}[D].$ Consider $W = W_{\operatorname{Exc}}(X)$. By \cite[Theorem~6.18(3),(5) and Lemma~6.23(1),(3)]{Markman_11}, the group $W$ is normal in a finite-index subgroup $\Gamma\subseteq O^+(\NS(X))$. Since the fixed subspace of $R_D$ is exactly $[D]^\perp$, we have $\operatorname{Fix}(W)=\bigcap_D[D]^\perp=R^\perp.$ Then the normality of $W$ shows that $R^\perp$ is invariant under $\Gamma$.

    Since the quadratic space $\NS(X)_\RR$ has signature $(1,\rho(X)-1)$ and $\rho(X)\geq 3$, by the same argument using Borel density in the proof of Theorem~\ref{thm:fibration_rank_3}, every $\Gamma$-invariant subspace is invariant under $\mathrm{SO}^+(\NS(X)_\RR)$. Since its standard representation is irreducible for $\rho(X)\geq 3$, we have that $R^\perp_\RR$ is either zero or all of $\NS(X)_\RR$. But since there is a prime exceptional divisor, the latter is impossible. Hence $R^\perp=0$. Then by the nondegeneracy of $q_X$, we have $R=\NS(X)_\QQ$.
\end{proof}

\begin{proposition}\label{prop:prime_divisor_dominate_B}
    Let $f:X\to B$ be a regular Lagrangian fibration on an irreducible hyperk\"ahler variety of dimension $2n$. If $X$ has a prime exceptional divisor, then it has one that dominates $B$.
\end{proposition}

\begin{proof}
    As before, we work over $\CC$. Let $L=f^*H$, where $H$ is an ample divisor on $B$, and let $A$ be an ample divisor on $X$. Since $\dim B=n$, we have $L^{n+1}=0$ for dimension reasons, and hence $L^{2n}=0$. By Fujiki's relation \cite[Lemma~2.5]{Kamenova_Lehn_25}, we have $\int_X L^{2n}=c_Xq_X(L)^n=0.$ Hence $q_X(L)=0$. Since $L$ is a nonzero class in the closure of the positive cone, we have $q_X(A,L)>0$. 

    Let $D$ be a prime divisor which does not dominate $B$. Since $f(D)$ is a proper closed subset of $B$, by the projection formula, we have 
    $$\int_X D L^n A^{n-1}=\int_B f_*(D A^{n-1})H^n=0.$$
    Moreover, by the polarized Fujiki relation \cite{Kamenova_Lehn_25}, we have 
    $$\binom{2n}{n}\int_X D L^n A^{n-1}=c_X2^n q_X(D,L)q_X(A,L)^{n-1}.$$
    Then since $q_X(A,L)>0$, we have $q_X(D,L)=0$.

    Assume that $X$ has a prime exceptional divisor $E$, but every such divisor satisfies $f(E)\subsetneq B$. If $\rho(X)\geq 3$, by Lemma~\ref{lem:maximal-rank}, the classes of prime exceptional divisors span $\NS(X)_\QQ$. By the computation above, $q_X(L,\NS(X)_\QQ)=0$, contradicting the nondegeneracy of the Beauville--Bogomolov form on $\NS(X)_\QQ$, since $L\neq 0$. If $\rho(X)=2$, then $L^\perp\cap\NS(X)_\QQ=\QQ L$, since $L$ is isotropic. Thus $L^\perp\cap\NS(X)_\QQ$ contains no negative-square classes, contradicting $q_X(E,L)=0$. Hence some prime exceptional divisor dominates $B$. 
\end{proof}

\begin{theorem}\label{thm:inductive_theorem}
    Let $f:X\to B$ be a regular Lagrangian fibration on an irreducible hyperk\"ahler variety of dimension $2n$ over $k$ with geometrically simple abelian fibre. Assume that rational points on each geometric prime exceptional divisor is potentially dense. If $\rho(X)\geq 3$ or if $\rho(X)  = 2, b_2(X)\geq 5$, then rational points on $X$ is potentially dense. 
\end{theorem}

\begin{proof}
    If $X_{\bar{k}}$ has no prime exceptional divisisors, then the result follows from Theorems~\ref{thm:fibration_rank_3} and~\ref{thm:fibration_rank_2}. Otherwise, by Proposition~\ref{prop:prime_divisor_dominate_B}, there is a prime exceptional divisor $D$ dominating $B$. After replacing $k$ with a finite extension, we may assume that $D$ is defined over $k$ and $D(k)$ is Zariski dense. Furthermore, after shrinking $B$, we may assume that $T=f^{-1}(B)$ is a torsor under an abelian scheme $A\to B$. By Lemma~\ref{lem:auxiliary-multiplication}, there is a finite \'etale surjective morphism $\varphi:T\to A$ over $B$. Let $M=D\cap T$. Then $\varphi|_M$ defines a section $s$ of $A_M\to M$.

    We show that $s$ is nondegenerate at the geometric generic point of $M$. Let $\eta$ be the generic point of $B$. Since $D$ dominates $B$, we have $\dim M_\eta=n-1>0$. Since $\varphi$ is finite, its restriction to $M_\eta$ has positive-dimensional image. Hence $s$ is not torsion at the geometric generic point of $M$, since otherwise its image would be contained in a finite torsion subgroup scheme of $A_\eta$. Since the geometric generic fibre of $A$ is simple, every nontorsion point is nondegenerate. Thus $s$ is nondegenerate at the geometric generic point of $M$.
    
    Since $M(k)$ is Zariski dense and $M\to B$ is dominant, applying Proposition~\ref{prop:parameter-core} to $T\to B$ with $M$ as the parameter space, we see that $T(k)$ is Zariski dense. Hence $X(k)$ is Zariski dense, proving that rational points on $X$ is potentially dense.
\end{proof}

\section{Conjecture ~\ref{conj:one-fibration} for each known deformation family}
In this section, we prove Conjecture~\ref{conj:one-fibration} for each known deformation family, namely hyperk\"ahler varieties of $K3^{[n]}$-type, generalized Kummer type, OG6 type, and OG10 type. The main idea is to construct a projective Shafarevich--Tate twist of the given variety which admits two distinct rational Lagrangian fibrations: one inherited from the original fibration and the other arising from a second movable isotropic ray. We then use a dominant rational map from the twist to the original variety to deduce potential density. We first recall the Shafarevich--Tate theory from \cite{Abasheva_26,huybrechts2025periodindexproblemhyperkahlermanifolds}.

\subsection{Shafarevich--Tate twists and lattices}
Let $f:X\to B$ be a regular Lagrangian fibration on a complex irreducible hyperk\"ahler variety. As explained in \cite[Definition~1.1.4]{Abasheva_26} (see also \cite[Section~1.1]{huybrechts2025periodindexproblemhyperkahlermanifolds}), we consider the identity-component subsheaf of the sheaf of relative automorphisms, denoted by $\Aut^0(X/B)$, and define its first cohomology group, called the Shafarevich--Tate group, by $\Sha(X/B) = H^1(B,\Aut^0(X/B)).$ We denote by $\Sha^0(X/B)$ the connected component of the identity in $\Sha(X/B)$. By the cocycle interpretation, each class $\alpha\in\Sha(X/B)$ determines a Shafarevich--Tate twist $X^\alpha\to B$ of the original fibration.

Furthermore, $\Sha^0(X/B)$ is closely related to the second cohomology of $X$. Consider the exact sequence of sheaves on $B$ 
$$0\to\Gamma\to f_*T_{X/B} \xrightarrow{\exp}\Aut^0(X/B)\to 0,$$
where $\Gamma$ is the kernel of the exponential map. Taking the long exact sequence in cohomology, we obtain
$$\Sha^0(X/B)\simeq H^1(B,f_*T_{X/B})/ \operatorname{im}\bigl(H^1(B,\Gamma)\bigr).$$
By \cite[Theorem~2.2.7]{Abasheva_26}, the holomorphic symplectic form and the Leray spectral sequence give an isogeny
$$\Sha^0(X/B)\sim H^{0,2}(X)/\pr^{0,2}(H^2(X,\ZZ)),$$
where $\pr^{0,2}$ denotes the Hodge projection.

\begin{remark}
    We recall the main results of \cite{Abasheva_26} that will be frequently used in our proofs. In particular, \cite[Theorem~2.2.10]{Abasheva_26} shows that the Shafarevich--Tate twists corresponding to classes in $\Sha^0(X/B)$ form a degenerate twistor family. By \cite[Theorem~A]{Abasheva_26}, these twists are K\"ahler. Moreover, when $X$ is projective, as we assume throughout this section, the twist $X^\alpha$ is projective if and only if $\alpha$ is torsion. Finally, the Lagrangian fibration on $X^\alpha$ determines a primitive integral class $e\in\NS(X^\alpha)$ generating its isotropic fibration ray \cite[Theorem 3.2.3]{Abasheva_26}.
\end{remark}

We next explain the relation between these Shafarevich--Tate classes and lattices. Let $\Lambda=H^2(X,\ZZ)$, equipped with the Beauville--Bogomolov form, and let $\Omega\in H^{2,0}(X)$ be a nonzero holomorphic symplectic form, normalized so that $(\Omega,\overline{\Omega})=1$. Then the Hodge projection is given by $\pr^{0,2}(a)=(a,\Omega)\overline{\Omega}$. Let $H$ be an ample line bundle on $B$ and put $\eta=f^*c_1(H)$. We denote by $e\in\NS(X)$ the primitive integral generator of the isotropic fibration ray. Thus $\eta=ce$ for some positive integer $c$.

By \cite[Theorem~2.2.10 and Section~5.1.5]{Abasheva_26}, the connected Shafarevich--Tate twists form a degenerate twistor family with marked periods $[\Omega+s\eta]$, where the parameter $s$ corresponds to $s\overline{\Omega}\in H^{0,2}(X)$ under the identification of \cite[Section~2.2.9]{Abasheva_26}. Writing $t=cs$, these periods become $[\Omega+te]$. We use this rescaled coordinate below and denote the corresponding class in $\Sha^0(X/B)$ by $[t]$. By \cite[Lemma~3.2.2]{Abasheva_26}, this class is torsion whenever $t/c=(a,\Omega)$ for some $a\in H^2(X,\QQ)$. In particular, $[(a,\Omega)]$ is torsion for every $a\in H^2(X,\QQ)$, since $(a,\Omega)/c=(a/c,\Omega)$.

We now turn to the study of integral lattices. We first present a lemma giving coordinates for a lattice and its localizations at all but finitely many primes.

\begin{lemma}
    Let $\Lambda$ be a nondegenerate integral lattice and let $e\in\Lambda$ be a nonzero primitive isotropic vector. Then there exists a rational isotropic vector $f_0\in\Lambda_\QQ$ such that $(e,f_0)=1$.

    Set $T_\QQ=\langle e,f_0\rangle^\perp\subset\Lambda_\QQ$ and $M=\Lambda\cap T_\QQ$. Then $\Lambda_\QQ=\QQ e\oplus\QQ f_0\oplus T_\QQ.$ Furthermore, for all but finitely many primes $p$, we have
    \begin{equation}\label{eq:prime-local-splitting}
    \Lambda\otimes\mathbb Z_{(p)} =\mathbb Z_{(p)}e \oplus\mathbb Z_{(p)}f_0\oplus(M\otimes\mathbb Z_{(p)}).
    \end{equation}
\end{lemma}

\begin{proof}
    Let $g\in\Lambda_\QQ$ be such that $(e,g)=1$. Define $f_0=g-\frac{q(g)}{2}e.$ By direct computation, using $q(e)=0$ and $(e,g)=1$, we have $q(f_0)=0$ and $(e,f_0)=1$.

    Next, choose a basis of $M$. Together with $e$ and $f_0$, it forms a basis of $\Lambda_\QQ$, since $M$ spans $T_\QQ$ over $\QQ$. Since this basis and an integral basis of $\Lambda$ are related by an invertible change-of-basis matrix over $\QQ$, the matrix and its inverse have entries in $\ZZ_{(p)}$ for all but finitely many primes $p$. Hence we obtain \eqref{eq:prime-local-splitting}.
\end{proof}

We apply the lemma to $\Lambda=H^2(X,\ZZ)$, equipped with the Beauville--Bogomolov form. Since $e$ is algebraic, we have $(\Omega,e)=0$. Writing $\Omega=\lambda e+\mu f_0+\sigma$ with $\sigma\in T_{\CC}$ and pairing with $e$, we obtain $\mu=0$. Thus $\Omega=\sigma+\lambda e$. Now choose a rational ample class $h=ue+rf_0+w$, where $w\in T_\QQ$. Since $e$ lies on the nef isotropic ray, we have $r=(h,e)>0$. Furthermore, since $(h,\Omega)=0$, we obtain $0=(w,\sigma)+r\lambda$. Setting $b=-\frac{w}{r}\in T_\QQ$, we have $\Omega=\sigma+(b,\sigma)e$. 

We define $N=\{z\in T_\QQ:(z,\sigma)=0\}$. Under the identification $T_\QQ\simeq e^\perp_\QQ/\QQ e$, we have 
$$N\simeq \bigl(\NS(X)_\QQ\cap e^\perp_\QQ\bigr)/\QQ e.$$
Since $\NS(X)_\QQ$ has signature $(1,\rho(X)-1)$ and $e$ is isotropic, the Hodge index theorem implies that $N$ is negative definite.

Furthermore, since the Beauville--Bogomolov form on $\Lambda_\QQ$ has signature $(3,b_2(X)-3)$ and $\langle e,f_0\rangle$ is a hyperbolic plane of signature $(1,1)$, the orthogonal complement $T_\QQ$ has signature $(2,b_2(X)-4)$. Thus $T_\QQ=N\oplus N^\perp$, where the orthogonal complement is taken inside $T_\QQ$. Since $N$ is negative definite, $N^\perp$ is nondegenerate with two positive directions. In particular, $\dim N^\perp\geq 2$.

In the next lemma, we show that the classes represented by $(a,\sigma)$ vanish in $\Sha^0(X/B)$ for all $a$ in a finite-index sublattice of $M$. This allows us to construct torsion twists of order dividing any prime.

\begin{lemma}\label{lem:twist-parameter-kernel}
    There exists an integer $m>0$ such that $[(a,\sigma)]=0$ in $\Sha^0(X/B)$ for every $a\in mM$. In particular, for every prime $p$ and every $a\in mM$, we have $\left[\frac{(a,\sigma)}{p}\right]\in\Sha^0(X/B)[p].$
\end{lemma}

\begin{proof}
    Since $(a,e)=0$ for $a\in M$, we have $(a,\Omega)=(a,\sigma)$. Then by \cite[Lemma~3.2.2]{Abasheva_26}, the corresponding class in $\Sha^0(X/B)$ is torsion. Choose an integral basis $a_1,\ldots,a_s$ of $M$, and let $m>0$ be a common multiple of the orders of the corresponding torsion classes. By additivity, $(a,\sigma)$ represents zero in $\Sha^0(X/B)$ whenever $a\in mM$. Therefore, for every prime $p$ and $a\in mM$, we have 
    $$p\left[\frac{(a,\sigma)}{p}\right]=[(a,\sigma)]=0.$$
\end{proof}

We now consider the projective torsion twists obtained from Lemma~\ref{lem:twist-parameter-kernel} and describe their N\'eron--Severi lattices. This will allow us to distinguish two isotropic rays.

\begin{lemma}\label{lem:torsion-twist-picard}
    Let $p$ be a prime and $a\in mM$. Define $\gamma=\left[\frac{(a,\sigma)}{p}\right]\in\Sha^0(X/B)[p],$ and $Y = X^\gamma$. Then $Y$ is projective with marked period $\Omega_Y=\sigma+\left(b+\frac{a}{p},\sigma\right)e.$ Its N\'eron--Severi lattice satisfies
    $$\NS(Y)_\QQ =\QQ e\oplus\QQ\left(f_0-b-\frac{a}{p}\right)\oplus N.$$
    In particular, $\rho(Y)=\dim_\QQ N+2$.

    Suppose further that \eqref{eq:prime-local-splitting} holds at $p$, that $b\in M\otimes\ZZ_{(p)}$, and that the orthogonal projections onto $N$ and $N^\perp$ preserve $M\otimes\ZZ_{(p)}$. If $t=\pr_{N^\perp}(pb+a)$ has nonzero reduction modulo $p$, then $p\mid \divs_{\NS(Y)} (e)$. 
\end{lemma}

\begin{proof}
    By Lemma~\ref{lem:twist-parameter-kernel}, we have $p\gamma=0$. Therefore, $Y$ is projective by \cite[Theorem~3.2.3]{Abasheva_26}, and it has marked period $\Omega_Y=\sigma+\left(b+\frac{a}{p},\sigma\right)e$ by \cite[Theorem~2.2.10 and Section~5.1.5]{Abasheva_26}.
    
    For $x=ue+rf_0+w\in\Lambda$, since $T_\QQ=\langle e,f_0\rangle^\perp$ and $(e,f_0)=1$, we have
    $$\begin{aligned}
    (x,\Omega_Y) &=\left(ue+rf_0+w,\, \sigma+\left(b+\frac ap,\sigma\right)e\right) =(w,\sigma)+r\left(b+\frac ap,\sigma\right)=\left(w+r\left(b+\frac ap\right),\sigma\right).
    \end{aligned}$$
    Hence $x$ is algebraic if and only if $w+r(b+a/p)\in N$. Writing
    $$x=ue+r\left(f_0-b-\frac ap\right) +\left(w+r\left(b+\frac ap\right)\right),$$
    we have $\NS(Y)_\QQ=\QQ e\oplus\QQ\left(f_0-b-\frac ap\right)\oplus N.$ In particular, $\rho(Y)=\dim_\QQ N+2$.

    Now assume that $t=\operatorname{pr}_{N^\perp}(pb+a)$ has nonzero reduction modulo $p$. For any $x=ue+rf_0+w\in\NS(Y)$, projecting the algebraicity condition onto $N^\perp$ gives $\operatorname{pr}_{N^\perp}(w)=-\frac rp t.$ By our assumptions on the local coordinates and orthogonal projections, $\operatorname{pr}_{N^\perp}(w)$ is $p$-integral. Since $t$ has nonzero reduction modulo $p$, we obtain $p\mid r$. As $r=(e,x)$, this holds for every $x\in\NS(Y)$, and hence $p\mid\operatorname{div}_{\NS(Y)}(e).$
\end{proof}

\subsection{Constructing a second fibration on the twist}
In this subsection, we construct a second rational Lagrangian fibration on the twist $Y=X^\gamma\to B$ constructed above and prove the theorem. We do this by finding another movable isotropic ray distinct from the ray generated by $e$. In particular, we consider two cases. When $N=0$, we can choose the twist so that its N\'eron--Severi lattice has rank $2$ and no walls. Hence the two isotropic boundary rays of its positive cone give two distinct rational Lagrangian fibrations. When $N\neq 0$, we use the lattice results above to construct another isotropic ray, and then apply reflections in prime exceptional divisors to move it into the movable cone. By comparing their divisibilities in the N\'eron--Severi lattice, we show that this ray is distinct from the ray generated by $e$.

\begin{proposition}\label{prop:prime-torsion-second-ray}
    Let $X$ be a complex projective irreducible hyperk\"ahler variety with $b_2(X)\geq 5$, and let $f:X\to B$ be a regular Lagrangian fibration. Let $e$ be the primitive integral generator of its nef isotropic ray. For every positive integer $d$ and every finite set $S$ of primes, there exist arbitrarily large primes $p\notin S$ satisfying $p\equiv 1\pmod d$ and classes $\gamma\in\Sha^0(X/B)[p]$ such that the projective twist $Y=X^\gamma$ has a second rational movable isotropic ray distinct from $\mathbb R_{\geq 0}e$.
\end{proposition}

\begin{proof}
    We keep the notation from the previous subsection. Since we want primes outside $S$, we may enlarge $S$ to contain all primes dividing $2m$, the denominators of $b$ and $f_0$, and the discriminants and indices needed for the decompositions \eqref{eq:prime-local-splitting} and $T_\QQ=N\oplus N^\perp$. Then for every prime $p\notin S$, the relevant lattices and orthogonal projections are defined over $\ZZ_{(p)}$, and the reductions of the corresponding quadratic forms are nondegenerate.

    Now assume that $N=0$. Since $\dim_{\FF_p}(M/pM)=b_2(X)-2\geq3$, the Chevalley--Warning theorem implies that $M/pM$ contains a nonzero isotropic vector. Since $p\nmid m$, the reduction of $mM$ surjects onto $M/pM$. We may therefore choose $a\in mM$ such that $t=pb+a\in M\otimes\ZZ_{(p)}$ satisfies $\bar t\neq0$ and $q(t)\equiv0\pmod p$.

    Let $\gamma=\left[\frac{(a,\sigma)}{p}\right]$ and $Y=X^\gamma$. By Lemma~\ref{lem:torsion-twist-picard}, we have $\rho(Y)=2$ and $p\mid r$ for every $x=ue+rf_0+w\in\NS(Y)$. Moreover, since $N=0$, the same lemma gives $w=-rt/p$. Writing $r=pR$ for some $R\in\ZZ$, we have
    $$q(x) = 2upR+R^2q(t)\in p\ZZ_{(p)}.$$
    Since $q(x)\in\ZZ$, we obtain $q(\NS(Y))\subset p\ZZ$. 

    By \cite[Theorem~3.17]{Amerik_Verbitsky_20}, the absolute values of the squares of primitive MBM classes are bounded by a constant depending only on the deformation type. We may choose $p$ larger than this bound. Then, since $q(\NS(Y))\subset p\ZZ$, every algebraic class of negative square has absolute square at least $p$. Hence $Y$ has no algebraic MBM classes. By \cite[Theorem~6.2]{Amerik_Verbitsky_15}, the nef cone of $Y$ is the closure of its positive cone. Since $\rho(Y)=2$ and $\NS(Y)$ contains the isotropic class $e$, the two isotropic boundary rays of its positive cone are rational and nef. Hence the other ray gives a second movable isotropic ray distinct from $\mathbb R_{\geq 0}e$.

    Now assume that $N\neq 0$. Choose $0\neq z_0\in M\cap N$, and enlarge $S$ to contain the primes dividing $q(z_0)$, which is nonzero since $N$ is negative definite. Then $L := M\cap N^\perp/p M\cap N^\perp$ is a nondegenerate quadratic space over $\FF_p$ of dimension at least $2$. Since $p$ is odd, such a quadratic form represents every nonzero scalar in $\FF_p$. 

    Since $p\nmid m$ and the orthogonal projections are integral over $\ZZ_{(p)}$, the reduction of $\pr_{N^\perp}(mM)$ surjects onto $L$. Hence we can choose $a\in mM$ such that $t = \pr_{N^\perp}(pb+a)$ satisfies $q(t)\equiv-q(z_0)\pmod p$. In particular, $q(t)$ is a unit in $\ZZ_{(p)}$. 

    Let $Y$ be as above. Consider the rational vector
    \begin{equation}\label{eq:second-null-vector}
    v=pf_0-t+z_0-\frac{q(t)+q(z_0)}{2p}e.
    \end{equation}
    Since $q(t)+q(z_0)\equiv 0\pmod p$ and $p$ is odd, the coefficient of $e$ lies in $\ZZ_{(p)}$. Hence $v\in\Lambda\otimes\ZZ_{(p)}$. Moreover, since $t\in N^\perp$ and $z_0\in N$, we have $(t,z_0)=0$. Thus, by direct computation,
    $$q(v)=q(t)+q(z_0)-2p\frac{q(t)+q(z_0)}{2p}=0.$$
    Furthermore, by Lemma~\ref{lem:torsion-twist-picard}, $v$ is algebraic, since
    $$(-t+z_0)+(pb+a)=\pr_N(pb+a)+z_0\in N.$$
    Hence $v\in\NS(Y)_\QQ$ is a nonzero isotropic vector.
    
    Furthermore, we have $(e,v)=p>0$ and $(v,z_0)=q(z_0)\in\ZZ_{(p)}^\times$. Since $v\in\Lambda\otimes\ZZ_{(p)}$, we can clear its denominators by multiplying it by a positive integer $c$ prime to $p$. Let $v_0$ be the primitive integral generator of this ray. Writing $cv=\ell v_0$ for some positive integer $\ell$, we have $c q(z_0)=\ell(v_0,z_0).$ Since $cq(z_0)$ is a $p$-adic unit, both $\ell$ and $(v_0,z_0)$ are $p$-adic units. As $z_0\in\NS(Y)$, we obtain $p\nmid\operatorname{div}_{\NS(Y)}(v_0).$ Furthermore, since $q(t)$ is a $p$-adic unit, $t$ has nonzero reduction modulo $p$. Hence, by Lemma~\ref{lem:torsion-twist-picard}, we have $p\mid\operatorname{div}_{\NS(Y)}(e).$

    Since $(e,v_0)>0$, $v_0$ lies in the closure of the same positive cone as $e$. By \cite[Proposition~3.4]{Riess_23}, there exists a composition of reflections in prime exceptional divisors sending $v_0$ to a class $v_1$ in the closed birational K\"ahler cone. These reflections are integral Hodge isometries by \cite[Proposition~3.3]{Riess_23}, and hence preserve primitivity and divisibility in $\NS(Y)$. Furthermore, by Lemma ~\ref{lem:torsion-twist-picard}, $v_1$ generates a ray distinct from that of $e$, and by \cite[Corollary~2.6]{Riess_23}, it lies in the closed movable cone. Finally, by Dirichlet's theorem, we can choose arbitrarily large primes $p\equiv 1\pmod d$ outside the finite set $S$. This completes the proof.
\end{proof}

We now construct a dominant rational map from the twist to $X$, using the multiplication-by-$p$ map on the generic abelian torsor, in the spirit of Lemma~\ref{lem:auxiliary-multiplication}.

\begin{lemma}\label{lem:torsion-twist-forward}
    Let $f:X\to B$ be a regular Lagrangian fibration of relative dimension $n$ on a complex projective irreducible hyperk\"ahler variety $X$. Let $d>0$ be the degree of the finite \'etale morphism constructed in Lemma~\ref{lem:auxiliary-multiplication} over a nonempty smooth open subset $U\subset B$. Assume that $r>0$ satisfies $r\equiv 1\pmod d$, and let $\gamma\in\Sha(X/B)[r]$ be such that the twist $Y:=X^\gamma$ is projective. Then there exists a dominant algebraic rational map $F: Y\dashrightarrow X$ over $B$, of generic degree $r^{2n}$, whose restriction to the smooth geometric generic fibre is multiplication by $r$ followed by a translation, after choosing origins. 
\end{lemma}

\begin{proof}
    Let $U\subset B$ be a nonempty open subset such that $T=f^{-1}(U)$ is a torsor under an abelian scheme $A\to U$. By Lemma~\ref{lem:auxiliary-multiplication}, there is a finite \'etale morphism $\varphi:T\to A$ such that $\varphi(t+u)=\varphi(t)+[\delta]u$, where $\delta>0$ is the degree of the multisection used in its construction. In particular, $\deg\varphi=\delta^{2n}=d$. Since $r\equiv1\pmod d$, we can write $r=1+m\delta$ for some integer $m$. Define $F_0:T\to T$ by $F_0(t)=t+[m]\varphi(t)$. Then for every $u\in A$, we have
    $$F_0(t+u)=t+u+[m]\bigl(\varphi(t)+[\delta]u\bigr)=F_0(t)+[r]u.$$
    Thus $F_0$ is fibrewise multiplication by $r$ followed by translation. Hence it is finite \'etale of degree $r^{2n}$ over $U$ and extends to a rational map $F_0:X\dashrightarrow X$ over $B$.

    By \cite[Lemma~3.1.2]{Abasheva_26}, we can represent $\gamma$ on an analytic open cover $\{U_i\}$ of $B$ by a cocycle of vertical automorphisms $\gamma_{ij}$ satisfying $r\gamma_{ij}=0$. Over the smooth locus, these automorphisms act by translations. Hence $F_0\circ\gamma_{ij}=F_0$, since $F_0$ is equivariant with respect to multiplication by $r$. As both sides are meromorphic, this equality extends to the overlaps. Therefore, the local maps $F_0$ glue to a meromorphic map $F:Y\dashrightarrow X$ over $B$.

    Since $X$ and $Y$ are projective, the graph of $F$ is algebraic by Chow's theorem. Hence $F$ is a dominant algebraic rational map, whose generic degree is $r^{2n}$.
\end{proof}

\subsection{Proof of the main theorem and applications}
In this subsection, we prove the main theorem and derive some immediate applications.

\begin{theorem}\label{thm:main}
    Let $X$ be an irreducible hyperk\"ahler variety over a number field $k$ of $K3^{[n]}$-type or generalized Kummer type with $n\geq 2$, or of OG6 or OG10 type. If, after a finite extension of $k$, $X$ admits a rational Lagrangian fibration, then rational points on $X$ is potentially dense.
\end{theorem}

\begin{proof}
    After replacing $k$ by a finite extension, we may assume that $X$ admits a rational Lagrangian fibration over $k$. By Lemma~\ref{lem:descent_finite_extension} and birational invariance of potential density, we may replace $X$ by a smooth projective birational model on which the Lagrangian fibration $f:X\to B$ is regular.

    The second Betti numbers of the four deformation types are respectively $23,7,8,24$ by \cite[Sections~6--7]{Beauville_83}, \cite{Rap07} and \cite{Rap08}. In particular, $b_2(X)\geq 5$. Let $\delta$ be the degree of a finite \'etale multisection supplied by Lemma~\ref{lem:auxiliary-multiplication}, and let $d=\delta^{2n}$ be the degree of the corresponding finite \'etale morphism to the acting abelian scheme. By Proposition~\ref{prop:prime-torsion-second-ray}, there is a prime $p\equiv 1\pmod d$ and a class $\gamma\in\Sha^0(X/B)[p]$ such that the projective twist $Y=X^\gamma$ has a rational movable isotropic ray distinct from the original ray generated by $e$. Since $Y$ is a connected Shafarevich--Tate twist, it has the same deformation type as $X$ by \cite[Theorem~2.2.10]{Abasheva_26}.

    Since $b_2(Y)\neq 4$, by \cite[Proposition~2.7]{Riess_23}, the second isotropic ray becomes nef on a smooth projective birational model $Y'$ of $Y$. Furthermore, by \cite[Theorem~4.5]{Riess_23}, \cite[Corollary~1.1]{Mat17}, \cite[Theorem~7.2]{MR21} and \cite[Theorem~2.2]{MO22}, respectively, this ray induces a rational Lagrangian fibration for each of the four deformation types. Since the original Lagrangian fibration is preserved by the twist, $Y$ admits two rational Lagrangian fibrations with distinct associated rays.

    By Lemma~\ref{lem:torsion-twist-forward}, there is a dominant algebraic rational map $F:Y\dashrightarrow X_\CC$ of generic degree $p^{2n}$. We now descend this construction to a finite extension of $k$. As in the proof of Lemma~\ref{lem:descent_finite_extension}, we can spread out $Y$, the smooth birational models carrying the two regular Lagrangian fibrations, their birational maps, and the dominant rational map $F$ to the fixed variety $X_{\bar{k}}$ over an integral scheme of finite type over $\bar{k}$. After shrinking the parameter scheme, we may assume that the fibres remain smooth projective hyperk\"ahler varieties, the fibrations remain Lagrangian with geometrically connected fibres, and the loci where the joint fibration map and $F$ have full differential rank remain nonempty. Choosing a closed point, we obtain a hyperk\"ahler variety $Z$ over a finite extension of $k$ admitting two rational Lagrangian fibrations with distinct associated rays and a dominant rational map to $X$. By Theorem~\ref{thm:two-fibrations}, rational points on $Z$ is potentially dense. Since it admits a dominant rational map to $X$, we conclude that rational points on $X$ is potentially dense.
\end{proof}

\begin{remark}
    Note that the one-fibration theorem holds for arbitrary deformation types as long as $b_2(X)\geq 5$. The only part of the proof that uses results specific to the four known deformation types is the realization of the second movable isotropic ray as a rational Lagrangian fibration, which follows in general from \cite[Theorem~1.1]{Engel_Mauri_26}.
\end{remark}

We now give some applications. 

\begin{corollary}\label{cor:Pic-5}
    Let $X$ be an irreducible hyperk\"ahler variety over a number field $k$ with $\rho(X_{\bar{k}})\geq 5$. Then rational points on $X$ is potentially dense.
\end{corollary}

\begin{proof}
    By the Hodge index theorem, $\NS(X_{\bar{k}})_\QQ$ has signature $(1,\rho(X_{\bar{k}})-1)$. Since $\rho(X_{\bar{k}})\geq 5$, by Meyer's theorem , it contains a nonzero rational isotropic class. After possibly changing its sign, by \cite[Proposition~3.4]{Riess_23}, we may apply reflections in prime exceptional divisors to move this class into the closed birational K\"ahler cone. Since $b_2(X)\geq \rho(X_{\bar{k}})+2\geq 7$, by \cite[Proposition~2.7]{Riess_23}, this class becomes nef on a smooth birational model. Then by the hyperk\"ahler SYZ theorem \cite[Theorem~1.1]{Engel_Mauri_26}, it defines a regular Lagrangian fibration. By Lemma~\ref{lem:descent_finite_extension}, we obtain a rational Lagrangian fibration after a finite extension of $k$.

    Finally, the proof of Theorem~\ref{thm:main} applies to arbitrary deformation types with $b_2(X)\geq 5$, using \cite[Theorem~1.1]{Engel_Mauri_26} to realize the second isotropic ray as a rational Lagrangian fibration. Hence rational points on $X$ is potentially dense.
\end{proof}

We next apply the theorem to the extended lattices of $K3^{[n]}$-type and generalized Kummer varieties. Let $H^2(X,\ZZ)\hookrightarrow\widetilde{\Lambda}_X$ denote the \emph{extended Hodge lattice}, with $H^2(X,\ZZ)^\perp=\ZZ v$ for some primitive algebraic vector $v$ of positive square. We use Markman's extended lattice in $K3^{[n]}$-type \cite[Section~9]{Markman_11} and Wieneck's extended lattice in generalized Kummer type \cite[Theorem~4.9]{Wie18}. Let $N_X=\widetilde{\Lambda}_X^{1,1}\cap\widetilde{\Lambda}_X$. Then we have the rational orthogonal decomposition $N_{X,\QQ}=\NS(X_{\bar{k}})_\QQ\oplus\QQ v.$ We define the \emph{Witt index} $\ind(N_{X,\QQ})$ to be the maximal dimension of a totally isotropic subspace of $N_{X,\QQ}$.

\begin{corollary}\label{thm:witt-two-density}
    Let $X/k$ be a projective irreducible hyperk\"ahler variety of $K3^{[n]}$-type or generalized Kummer type, with $n\geq 2$. If $\ind(N_{X,\QQ})=2$, then rational points on $X$ is potentially dense.
\end{corollary}

\begin{proof}
    Choose a two-dimensional totally isotropic subspace $L\subset N_{X,\QQ}$. Since $v^2>0$, the restriction of $(v,-)$ to $L$ has a nonzero kernel. Thus $L\cap v^\perp$ contains a nonzero isotropic vector. By the rational orthogonal decomposition, this vector lies in $\NS(X_{\bar{k}})_\QQ$. After possibly changing its sign, by \cite[Proposition~3.4]{Riess_23}, we can move this class into the closed birational K\"ahler cone using reflections in prime exceptional divisors. By \cite[Theorem~4.5]{Riess_23} and \cite[Corollary~1.1]{Mat17}, the resulting isotropic ray induces a rational Lagrangian fibration on a smooth birational model of $X_\CC$. By Lemma~\ref{lem:descent_finite_extension}, this fibration can be realized after a finite extension of $k$. Potential density then follows from the one-fibration theorem proved above.
\end{proof}

\appendix
\section{An example in arbitrary dimensions}
In this appendix, we construct examples of hyperk\"ahler varieties admitting one Lagrangian fibration in arbitrary dimensions satifying the potential density property, using dominant rational maps from products of abelian varieties and K3 surfaces. 

We first show that the tangent plane of a K3 surface admitting no genus-one fibration is spanned by the tangent lines of a covering family of genus-one curves.

\begin{lemma}\label{lem:elliptic_web_directions}
    Let $S$ be a projective K3 surface over an algebraically closed field of characteristic zero admitting no genus-one fibration. Let $\cE\to C$ be a covering family of smooth genus-one curves, where $C$ is an integral curve, whose general fibre maps birationally onto an integral curve $E_a\subset S$. Then the tangent lines $T_xE_a$ of the curves $E_a$ passing through a geometric generic point $x\in S$ span $T_xS$.
\end{lemma}

\begin{proof}
    We may replace $C$ by the normalization of its image in $\operatorname{Chow}_{1,d}(S)$, so that distinct general parameters correspond to distinct curves in the Chow variety. Then we normalize the universal incidence and restrict to a smooth family. Since the family covers $S$ and both $\cE$ and $S$ have dimension two, the evaluation map $e:\cE\to S$ is dominant and generically finite. After compactifying the family and shrinking to a nonempty open subset $U\subset S$, the evaluation map $e$ is finite \'etale over $U$, and the image curves are smooth at their points in $U$.

    Assume that the tangent lines $T_xE_a$ through a geometric generic point $x$ do not span $T_xS$. Then they coincide and, after shrinking $U$, define a rank-one subbundle of $T_U$. Every curve $E_a$ meeting $U$ is tangent to this distribution at its general points. Since a rank-one distribution is integrable, by uniqueness of integral curves in characteristic zero, any two members through a general point agree as germs, hence as integral curves. Since distinct general parameters correspond to distinct curves, the evaluation map is generically one-to-one. Hence the inverse of the evaluation map, followed by the projection to $C$, gives a rational genus-one fibration on $S$.

    Resolving this rational map and taking its Stein factorization, we obtain a morphism $h:\widetilde S\to C'$ with connected fibres, where $\widetilde S\to S$ is a composition of blowups and $C'$ is a smooth projective curve. Its general fibre $F$ is a smooth genus-one curve, so $F^2=0$ and, by adjunction, $K_{\widetilde S}\cdot F=0$. Since $S$ is a K3 surface, we have $K_{\widetilde S}=\sum_i a_iE_i$, where the $E_i$ are the exceptional curves and $a_i>0$. As $F$ is nef, we have $E_i\cdot F\geq 0$ for every $i$. Hence $E_i\cdot F=0$ for every exceptional curve $E_i$, so all exceptional curves are vertical for $h$. Therefore, $h$ descends to a genus-one fibration on $S$, contradicting our assumption.
\end{proof}

Next, we give a lemma about choosing a basis of the quotient by a Lagrangian subspace from an orthogonal direct sum of symplectic planes. This will be used when we consider the tangent space of the product.

\begin{lemma}\label{lem:product-direction}
    Let $W=\bigoplus_{i=1}^n W_i$ be an orthogonal direct sum of symplectic planes, and let $K\subset W$ be a Lagrangian subspace. For each $i$, choose a finite spanning set $P_i\subset W_i$. Then there exist $v_i\in P_i$ whose images in $W/K$ form a basis of $W/K$.
\end{lemma}

\begin{proof}
    Let $\pi:W\to W/K$ be the quotient map. For any $I\subset\{1,\ldots,n\}$, set $W_I=\bigoplus_{i\in I}W_i$. Since $K$ is Lagrangian, $K\cap W_I$ is an isotropic subspace of $W_I$. Hence $\dim(K\cap W_I)\leq |I|$. Since each $P_i$ spans $W_i$, we have
    $$\dim\left\langle\pi(P_i):i\in I\right\rangle =\dim\pi(W_I) =2|I|-\dim(K\cap W_I)\geq |I|.$$
    Then by Rado's theorem \cite{Rado1942}, there exist $v_i\in P_i$ whose images $\pi(v_i)$ are linearly independent. Since $\dim(W/K)=n$, these images form a basis.
\end{proof}

Next, we show that a product of abelian varieties and K3 surfaces admits a dominant generically finite evaluation map from a family of products of genus-one curves, abelian varieties, and K3 surfaces admitting genus-one fibrations.

\begin{lemma}\label{lem:dominant_from_genus_one}
    Let $P=A\times S_1\times\cdots\times S_r$, where $A$ is an abelian variety of dimension $2(n-r)$ and the $S_i$ are projective K3 surfaces over an algebraically closed field of characteristic zero. Assume that $P$ admits a symplectic form
    $$\omega=\sigma_A+\sum_{i=1}^r c_i\operatorname{pr}_i^*\sigma_i,$$ 
    where $\sigma_A$ is symplectic and $c_i\neq 0$. Let $h:P\dashrightarrow B$ be a dominant rational map to a smooth variety of dimension $n$, whose geometric generic tangent kernel is Lagrangian. Then there exists a family $\mathcal M\to\Gamma$ of products of genus-one curves, abelian varieties, and K3 surfaces admitting genus-one fibrations, with a dominant generically finite evaluation map $e:\mathcal M\to P$, such that the image under $e$ of the locus where the restriction of $h\circ e$ to a member has differential of rank $n$ contains a nonempty open subset of $P$.
\end{lemma}

\begin{proof}
    On each K3 factor not admitting a genus-one fibration, choose a covering family of smooth genus-one curves as in \cite[Theorem~4.1]{Hassett_Tschinkel_00}. By Lemma~\ref{lem:elliptic_web_directions}, the tangent lines of these curves through a geometric generic point span the tangent plane of each factor.

    Let $p\in P$ be a geometric generic point where $h$ is defined, and let $K=\ker dh_p$. We decompose the tangent space of $P$ into an orthogonal sum of symplectic planes, taking the tangent planes of the K3 factors and a symplectic decomposition of the tangent space of $A$. Hence $T_pP=\bigoplus_{i=1}^n W_i.$ For each nonelliptic K3 factor, we choose a finite spanning set of tangent vectors from the covering family, and for the remaining planes we choose any finite spanning set. By Lemma~\ref{lem:product-direction}, we can choose one vector in each $W_i$ such that their images form a basis of $T_pP/K$.

    Choosing the corresponding genus-one curves and keeping the abelian factor and the K3 factors admitting genus-one fibrations, we obtain a product $M$ whose tangent space at a point above $p$ contains these vectors. Hence the differential of $h|_M$ has rank $n$, so the rational map $M\dashrightarrow B$ is dominant.

    Consider the universal family $\mathcal M\to\Gamma$ of such products, with its evaluation map $e:\mathcal M\to P$. Since the covering families dominate their respective K3 factors generically finitely, $e$ is dominant and generically finite. The locus $\mathcal M^\circ$ where $h\circ e$ is defined and its differential along the fibres of $\mathcal M\to\Gamma$ has rank $n$ is a nonempty open subset of $\mathcal M$. Since this condition holds above a geometric generic point of $P$, the image $e(\mathcal M^\circ)$ is dense in $P$. By Chevalley's theorem \cite[Exercise~II.3.19]{Hartshorne_77}, it contains a nonempty open subset of $P$.
\end{proof}

Now we prove our main theorem using the above lemmas.

\begin{theorem}\label{thm:example}
    Let $X$ be an irreducible hyperk\"ahler variety of dimension $2n$ over a number field $k$ admitting a rational Lagrangian fibration. Assume that after a finite extension of $k$, there is a dominant generically finite rational map
    $$P=A\times S_1\times\cdots\times S_r\dashrightarrow X,$$ 
    where $A$ is an abelian variety of dimension $2(n-r)$ and the $S_i$ are projective K3 surfaces. Then rational points on $X$ is potentially dense.
\end{theorem}

\begin{proof}
    Since potential density is a birational invariant, we may replace $X$ by its birational model on which the Lagrangian fibration is regular. Assume that $f:X\to B$ is regular and let $g:P\dashrightarrow X$ be the dominant generically finite rational map. After a finite extension, we may assume that everything is defined over $k$. If every $S_i$ admits a genus-one fibration, then we are done by \cite{Bogomolov_Tschinkel_00} and \cite[Proposition~3.1]{Hassett_Tschinkel_00}, since rational points on $P$ is potentially dense. Thus we assume that at least one K3 factor has no genus-one fibration.

    Over $\bar{k}$, let $\sigma_X$ be the symplectic form on $X$. Then $g^*\sigma_X$ extends to a regular $2$-form on $P$ by Hartogs' theorem, since $P$ is smooth and the indeterminacy locus of $g$ has codimension at least two. By the K\"unneth formula and $H^0(S_i,\Omega_{S_i}^1)=0$, we have
    $$g^*\sigma_X=\sigma_A+\sum_{i=1}^r c_i\operatorname{pr}_i^*\sigma_i,$$
    where $\sigma_A$ is a translation-invariant $2$-form on $A$ and $\sigma_i$ is a symplectic form on $S_i$. Since $g$ is generically finite, $g^*\sigma_X$ is generically nondegenerate. Hence all $c_i$ are nonzero and $\sigma_A$ is nondegenerate. Therefore, $g^*\sigma_X$ is nondegenerate everywhere on $P$, and $dg$ is an isomorphism wherever $g$ is defined.

    As before, after shrinking $B$ to a nonempty smooth open subset $U$, we may assume that $T=f^{-1}(U)\to U$ is a torsor under an abelian scheme $A_0\to U$ by \cite[Lemma~1.16]{Huybrechts_Mauri_22} and \cite[Theorem~5.3]{Laurent_Schroer_24}. Then by Lemma~\ref{lem:auxiliary-multiplication}, there is a finite \'etale surjective $U$-morphism $\varphi:T\to A_0$. Let $p_0\in P$ be a geometric generic point where $g$ is defined and $g(p_0)\in T$. Since $dg$ is an isomorphism and the fibres of $f$ are Lagrangian, the subspace $K=\ker d(f\circ g)_{p_0}\subset T_{p_0}P$ is Lagrangian with respect to $g^*\sigma_X$.

    Applying Lemma~\ref{lem:dominant_from_genus_one} to $h=f\circ g:P\dashrightarrow U$, we obtain a family $\mathcal M\to\Gamma$ of products of genus-one curves, abelian varieties, and K3 surfaces admitting genus-one fibrations, with a dominant generically finite evaluation map $e:\mathcal M\to P$. Let $\mathcal M^\circ$ be the nonempty open subset where $g\circ e$ is defined, its image lies in $T$, and the differential of $f\circ g\circ e$ along the fibres of $\mathcal M\to\Gamma$ has rank $n$. By the lemma, the image of this locus under $e$ contains a nonempty open subset of $P$. Since $g$ is dominant, by Chevalley's theorem \cite[Exercise~II.3.19]{Hartshorne_77}, the image of $\mathcal M^\circ$ under $g\circ e$ contains a nonempty open subset $\Omega\subset T$.

    Choose an algebraic point $b\in U$ such that $\Omega_b$ is nonempty. Since $\varphi_b$ is finite \'etale and surjective, $\varphi_b(\Omega_b)$ is a nonempty open subset of $A_{0,b}$. By \cite[Proposition~3.1]{Hassett_Tschinkel_00}, after a finite extension, $A_{0,b}$ has a nondegenerate point $a$. Hence we can choose a nonzero multiple $q=[m]a\in\varphi_b(\Omega_b)$, which is still nondegenerate. After another finite extension, we can choose $x\in\Omega_b$ such that $\varphi_b(x)=q$.
    
    Then lifting $x$ to $\mathcal M^\circ$, we can fix the corresponding genus-one curves from the covering families. Let $M$ be the product of these genus-one curves, the K3 factors admitting genus-one fibrations, and the abelian factor $A$. Then the evaluation map induces a rational map $j:M\dashrightarrow T$ with a point $p\in M$ such that $j(p)=x$. Since the differential of $f\circ j$ has rank $n$ at $p$, the composition $f\circ j:M\dashrightarrow U$ is dominant.
    
    Rational points on each of these factors is potentially dense by \cite[Proposition~3.1]{Hassett_Tschinkel_00} and \cite{Bogomolov_Tschinkel_00}. Hence rational points on $M$ is potentially dense after a finite extension. The composition $\varphi\circ j$ induces a rational section $s:M\dashrightarrow A_{0,M}$. At the chosen lift $p$ of $x$, the section has value $s(p)=\varphi_b(x)$, which is nondegenerate. Hence by Lemma~\ref{lem:cyclic-specialization}, $s$ is nondegenerate at the geometric generic point of $M$. Then by Proposition~\ref{prop:parameter-core}, rational points on $X$ is potentially dense. 
\end{proof}

\printbibliography

\end{document}